\documentclass[a4paper,12pt]{article}

\usepackage{amsmath,amsfonts,amssymb,amsthm}
\theoremstyle{plain}
\newtheorem{thm}{Theorem}
\newtheorem*{thm*}{Theorem}

\newtheorem*{lem*}{Lemma}

\theoremstyle{definition}
\newtheorem*{example}{Example}
\newtheorem*{Def*}{Definition}

\newcommand{\reftit}{\textit}    % reference title
\newcommand{\refis}{\textbf}     % reference issue no.

\begin{document}

\title{Exclusion Processes with Multiple Interactions.}
\author{Yevgeniy Kovchegov \footnote{This research was supported in part by NSF VIGRE Grant DMS
 9983726 at UCLA}}
\date{ }
\maketitle

\begin{abstract}
             We introduce the mathematical theory of the particle systems that interact
               via permutations, where the transition rates are assigned not to the jumps from
               a site to a site,
               but to the permutations themselves. This permutation processes can be viewed 
              as a generalization of the symmetric exclusion processes, where 
              particles interact via transpositions. 
              The duality and coupling techniques for the processes are described, the needed
              conditions for them to apply are established.  
              The stationary distributions of the permutation processes are explored for
              translation invariant cases.
\end{abstract}

\section{Introduction.}
 We begin by reformulating the general setup of the symmetric exclusion process.
 Let $S$ be a general countable set, and $p(x,y)$ be transition probabilities for
 a Markov chain on $S$. Let $\eta_t$ denote a continuous time Feller process with 
 values in $\{0,1\}^S$, where $\eta_t(x)=1$ when the site $x \in S$ is occupied
 by a particle at time $t$ while $\eta_t(x)=0$ means the site is empty at time $t$.
 The exclusion process is a fine example of a Markovian interacting particle system, with the
 name justified by the transition rates 
 $$\eta \rightarrow \eta_{x,y} \text{\quad at rate \quad} p(x,y) \text{\quad if \quad}
    \eta(x)=1, \eta(y)=0,$$
 where for $\eta \in \{0,1\}^S$, $\eta_{x,y}(u) := \eta(u)$ when $u \not\in \{x,y\}$,
 $\eta_{x,y}(x) := \eta(y)$ and $\eta_{x,y}(y) := \eta(x)$.
 The condition
 $$\sup_{y \in S} \sum_x p(x,y) < \infty$$
 is sufficient to guarantee that  the exclusion process $\eta_t$ is indeed a well defined
 Feller process. We refer the reader to \cite{liggett1} and \cite{liggett2} for a complete
 rigorous treatment of the subject.

 The exclusion process is {\it symmetric} if $p(y,x)=p(x,y)$ for all $x,y \in S$.
 In this case we can reformulate the process by considering all the transpositions
 $\tau_{x,y}$. For each transposition $\tau_{x,y}$ ($x,y \in S$, $x \not= y$) we will assign 
 the corresponding rate $q(\tau_{x,y})=p(y,x)=p(x,y)$ at which the transposition occurs:
 $$\eta \rightarrow \tau_{x,y}(\eta) \text{\quad at rate \quad} q(\tau_{x,y}),$$
 where $\tau_{x,y}(\eta) := \eta_{x,y}$.
 It was suggested  to the author by Tom Liggett to study the natural generalization
 of the process that arises with the above reformulation. 
 Liggett's idea was to assign the rates not to the particles
 inhabiting the space $S$, but to the various permutations of finitely many points of $S$. 
 Namely, we can consider other permutations besides the transpositions.
 We let $\Sigma$ be the set of all such permutations with {\bf positive} rate.
 If $\sigma \in \Sigma$, we let
 $$Range(\sigma)=\{x \in S \text{ :\quad} \sigma(x) \not= x\}.$$
 For each $\eta \in \{0,1\}^S$, let $\sigma(\eta)$ be the new configuration of particles
 after the permutation $\sigma$ was applied to $\eta$:
 $$\sigma(\eta)(x)=\eta(\sigma^{-1}(x)) \text{\quad for all \quad} x \in S.$$
 Observe that we only permute the particles inside $Range(\sigma)$.
 
 Now, we want to construct a continuous time Feller process, where rates $q(\sigma)$
 ($\sigma \in \Sigma$) are assigned so that
 $$\eta \rightarrow \sigma(\eta) \text{\quad at rate \quad} q(\sigma).$$ 
 
 \begin{example}
 Let $S=\mathbb{Z}$, and $\Sigma=\bigcup_{x \in \mathbb{Z}} 
    \big\{  \sigma_x :=  \overline{(x,x+1,x+2)}, \sigma_x^2=\sigma_x^{-1} \big\}$
 consists of all the three-cycles of consecutive integers. As we will see later, the three-cycles are
 very special for the theory of ``permutation" processes described in this manuscript.
 \end{example}
 
 First we would like to mention some of the results from
 the theory of exclusion processes 
 that we will extend to the newly introduced permutation processes.
 For consistency we will use the notations of \cite{liggett1} and \cite{liggett2}. We let  $\mathcal{I}$
 denote the class of stationary distributions for the given Feller process. As the set $\mathcal{I}$
 is convex, we will denote by $\mathcal{I}_e$ the set of all the extreme points of $\mathcal{I}$.  
 The results that we want to generalize are the two theorems given below.
 Consider the case of  $S=\mathbb{Z}^d$ with shift-invariant random walk rates (e.g. $p(x,y)=p(0,y-x)$). 
 The first theorem was proved 
 by F.Spitzer (see \cite{spitzer74}) in the recurrent case and by T.Liggett
 in the transient case (see \cite{liggett3}).
 \begin{thm*}
 For the symmetric exclusion process, $\mathcal{I}_e=\{ \nu_{\rho}: 0 \leq \rho \leq 1 \}$,
 where $\nu_{\rho}$ is the homogeneous product measure on $\{0,1\}^S$ with marginal probability   
 $\rho$ (e.g. \\$\nu_{\rho}\{ \eta: \eta := 1 \text{ on } A\}=\rho^{|A|}$ for any $A \subset S$).
 \end{thm*}
 Let  $\mathcal{S}$ denote the class of the shift invariant probability measures on $\{0,1\}^S$,
 and $(\mathcal{I} \cap \mathcal{S})_e$ the set of all extreme points of
 $(\mathcal{I} \cap \mathcal{S})$. Next theorem was proved in \cite{liggett4} by T.Liggett.
 A special case of it was proved by R.Holley in \cite{holley}.
 \begin{thm*}
 For the general exclusion process, $(\mathcal{I} \cap \mathcal{S})_e=\{ \nu_{\rho}: 0 \leq \rho \leq 1\}$.
 \end{thm*}

 As it was the case with the exclusion processes, coupling method will play the crucial
 role in proving the analogues of the above results for the permutation processes. The difficult
 part was to construct the right types of couplings for the corresponding proof to work.

\subsection{Existence of the process. The permutation law.} 
 We need to formalize the construction of the permutation process. For a configuration $\eta \in \{0,1\}^S$ and 
 a permutation $\sigma \in \Sigma$, $\sigma(\eta)$ defined as
 $$\sigma(\eta)(x) := \eta(\sigma^{-1}(x)) \quad \text{for all } x \in S$$
  is the resulting configuration after the permutation $\sigma$ is applied.
  For any cylinder function $f$ (i.e. a function $f(\eta)$ from $\{0,1\}^S$ to $\mathbb{R}$ that depends
  on finitely many sites in $S$), let
  $$\Omega f(\eta) := \sum_{\sigma \in \Sigma} q(\sigma) [f(\sigma(\eta))-f(\eta)].$$

  Now, we have to guarantee that the permutation process $\eta_t$ with generator $\Omega$
  is a well defined Feller process. For this, by Theorem 3.9 of Chapter I in \cite{liggett1}
  (see also the conditions (3.3) and (3.8) there),
  it is sufficient to assume that the rates $q(\sigma)$ are such that for every $x \in S$,
  \begin{eqnarray} \label{PL}
  M_{PL} := \sup_{x \in S} \sum_{\sigma : x \in Range(\sigma)} q(\sigma) <\infty.
  \end{eqnarray}
  Then the semigroup $S_t$ of the permutation process $\eta_t$, generated by such $\Omega$,
  is well defined. Such process will then be said to obey the permutation law (\ref{PL}).
  
  Throughout the paper we require that the random walk generated by the permutations 
  $\{ \sigma \in \Sigma \}$ is \textbf{irreducible}.
  That is that for every $x$ and $y$ in $S$ there is a sequence
  $\sigma_1,\dots,\sigma_k \in \Sigma$ with $q(\sigma_i)>0$ for all $i=1,\dots,k$
  such that $\sigma_k\circ ...\circ \sigma_1(x)=y$.

\subsection{Duality.}
  For a nonnegative continuous function $H(\eta, \zeta)$ of two variables, the Markov
  processes $\eta_t$ and $\zeta_t$ are said to be dual with respect to $H(\cdot,\cdot)$
  if $$E^{\eta} H(\eta_t, \zeta)=E^{\zeta}H(\eta, \zeta_t)$$
  for all $\eta, \zeta$ and all $t \geq 0$.
  \\
  Let for a configuration $\eta \in \{0,1\}^S$ and a set $A \subset S$,
  $$H(\eta, A)=\prod_{x \in A} \eta(x).$$

  \begin{eqnarray*}
  \Omega H(\cdot, A)(\eta) 
  & = & \sum_{\sigma \in \Sigma} q(\sigma) [H(\sigma(\eta),A)-H(\eta, A)] \\
  & = & \sum_{\sigma \in \Sigma} q(\sigma) [H(\eta, \sigma^{-1}(A))-H(\eta, A)] \\
  & = & \sum_{\sigma \in \Sigma} q(\sigma) [H(\eta, \sigma(A))-H(\eta, A)], 
  \end{eqnarray*}  
  where the last line is true whenever 
  \begin{eqnarray} \label{duality}
  q(\sigma)=q(\sigma^{-1}) \mbox{  for all  } \sigma \in \Sigma.
  \end{eqnarray}
  Then $$\Omega H(\cdot, A)(\eta)=\Omega H(\eta,\cdot)(A),$$ and the permutation processes
  $\eta_t$ and $A_t$ with $\eta_0=\eta$ and $A_0=A$ are dual with respect to $H$.
  So the permutation process satisfying (\ref{duality}) is {\it self-dual}. Therefore
    \begin{eqnarray} \label{dualEq}
    P^{\eta}[\eta_t := 1 \text{ on } A]=P^{A}[\eta := 1 \text{ on } A_t].
    \end{eqnarray}
  The condition (\ref{duality}) is essential in order to have a useful duality. From now on we will
  say that the Feller process is a {\it symmetric} permutation process whenever the above condition
  (\ref{duality}) is satisfied. Observe that in this case the process is analogous to the symmetric
  exclusion process, where the corresponding self-duality was indispensable and is the reason
  why the symmetric exclusion was so successfully studied (see Chapter VIII of \cite{liggett1} 
  and Part III of \cite{liggett2}). 

\section{Symmetric permutation processes.}

 For the rest of the paper we restrict ourselves to studying {\bf permutation processes
 on $S=\mathbb{Z}^d$ and the rates $q(\sigma)$ are assumed to be shift invariant. }
 We will also assume that the rates 
   $q(\sigma)$, for all $\sigma \in \Sigma$, satisfy the following two conditions.
   First
  \begin{eqnarray} \label{I}
  \sup_{\sigma \in \Sigma } \big| Range(\sigma) \big| <\infty ,
  \end{eqnarray}
  where $|\cdot |$ denotes the cardinality.
  Second, if $\sigma_2$ is a finite permutation of elements in $S$ such that
  $Range(\sigma_2)=Range(\sigma_1)$ for some $\sigma_1 \in \Sigma$, then
  $\sigma_2 \in \Sigma$, and
  \begin{eqnarray} \label{II}
  \sup_{\sigma_1, \sigma_2 \in \Sigma: Range(\sigma_1)=Range(\sigma_2)}
  \Big| {q(\sigma_1) \over q(\sigma_2)} \Big| <\infty .
  \end{eqnarray}
  From now on, we let $M_I$ denote the max in (\ref{I}) and $M_{II}$ denote the sup
  in (\ref{II}).
  It should be mentioned that the second condition (\ref{II}) is stricter than it needs to be.
  We only need $\Sigma$ to be the class of permutations where for the same range set,
  any ordering (word) of 1's and 0's on the range can be permuted into any other
  ordering with the same number of 1's and 0's by applying a permutation from that class.  
 In this section we assume that the process satisfies the duality conditions (\ref{duality}).
 We will prove 
 \begin{thm} \label{symmetric}
 For the symmetric permutation processes, $\mathcal{I}_e=\{ \nu_{\rho}: 0 \leq \rho \leq 1 \}$,
 where $\nu_{\rho}$ is the homogeneous product measure on $\{0,1\}^S$ with marginal probability   
 $\rho$ (e.g. \\$\nu_{\rho}\{ \eta: \eta := 1 \text{ on } A\}=\rho^{|A|}$).
 \end{thm}
 The notion of a bounded harmonic function for a Markov chain can be adapted to permutation
 processes. We will say that a bounded function $f:S \rightarrow \mathbb{R}$ is {\it harmonic}
 if for a permutation process $\eta_t$ and each $t>0$, 
 $f(\eta)=\sum_{\zeta \in \{0,1\}^S} P^{\eta}[\eta_t=\zeta] f(\zeta)$. We refer the reader to Chapter I
 of \cite{liggett1} for more on Markov processes, their semigroups and construction of interacting particle
 systems.
 We will need the following  
 \begin{thm} \label{harmonic}
 If $f$ is a bounded harmonic function for the well defined finite permutation process $A_t$,
 then $f$ is constant on $\{A: |A|=n\}$ for each given integer $n \geq 1$.
 \end{thm}
 As it was the case for symmetric exclusion, Theorem \ref{symmetric} follows from 
 Theorem \ref{harmonic} and the duality of the process (see \cite{liggett1}, Chapter VIII).
 The proof of Theorem \ref{symmetric} echos bit to bit the corresponding proof in case 
 of the symmetric exclusion process.
 However, we will briefly go through it. We assume that we already have Theorem \ref{harmonic}.
 
 {\it Proof of Theorem \ref{symmetric}:} A probability measure $\mu$ on $\{0,1\}^S$ is called
 {\it exchangeable} if for any finite $A \subset S$, $\mu\{\eta: \eta := 1 \text{ on }A\}$ 
 is a function of cardinality $|A|$ of
 $A$. By de Finetti's Theorem, if $S$ is infinite, then every exchangeable measure is
 a mixture of the homogeneous product measures $\nu_{\rho}$. Therefore
 Theorem \ref{symmetric} holds if and only if $\mathcal{I}$ agrees with the set of 
 exchangeable probability measures.
  
 The duality equation (\ref{dualEq}) implies
 \begin{eqnarray*}
 \mu S_t \{\eta: \eta := 1 \text{ on }A\} & = & \int P^{\eta} [\eta_t := 1 \text{ on } A] d\mu \\
 & = & \int P^A [\eta := 1 \text{ on } A_t] d\mu \\
 & = & \sum_B P^A[A_t=B] \mu\{\eta := 1 \text{ on } B \}.
 \end{eqnarray*}
 Thus every exchangeable measure is stationary. Now, if $\mu \in \mathcal{I}$, then 
 $\mu S_t=\mu$ (for all $t$), so by the above equation, $f(A)=\mu\{\eta := 1 \text{ on } A \}$
 is harmonic for $A_t$. Hence Theorem \ref{harmonic} implies that $\mu$ is exchangeable.
 $\square$
 
 The proof of Theorem \ref{harmonic} is different for the processes with recurrent and transient
 rates. We will do both.

\subsection{Recurrent case.}

  By \textbf{recurrence} here we mean the recurrence of $I_1(t) -I_2(t)$, where $I_1(t)$ and
  $I_2(t)$ are independent one-point processes moving according to the permutation law
  as described in the introduction. For the rest of the subsection we will assume that the process
  is recurrent. 
   
 As it was the case with the symmetric exclusion processes, to prove Theorem \ref{harmonic} 
 for the recurrent case, it is enough to construct a {\it successful} coupling of two copies
 $A_t$ and $B_t$ of the permutation process with initial states $A_0$ and $B_0$ of
 the same cardinality $n$ that coincide at all but two sites of $S$ (e.g.
 $|A_0 \cap B_0|=n-1$). By successful coupling, we mean
 $$P[A_t=B_t \text{  for all } t \text{ beyond some time }]=1.$$
 If $f$ is a bounded harmonic function for the finite permutation process, for which
 we can construct a successful coupled process (see above), then
 $$|f(A_0)-f(B_0)|=|Ef(A_t)-Ef(B_t)| \leq E|f(A_t)-f(B_t)| $$
 $$ \leq \|f\|P[A_t \not= B_t].$$
 Letting $t$ go to infinity, we get $f(A_0)=f(B_0)$, proving Theorem \ref{harmonic}
 for the case when there are only two discrepancies between $A_0$ and $B_0$
 (the cardinalities $|A_0|$=$|B_0|$, and $|A_0 \cap B_0|=|A_0|-1$).
 By induction, Theorem \ref{harmonic} holds for all $A_0$ and $B_0$ of the same
 cardinality.
 
 Now we need to construct a successful coupling with the property that
 $$P[A_t = B_t \mbox{ for all }t \mbox{ beyond some point }]=1.$$
 The points in ${\{(A_t \cup B_t) \backslash (A_t \cap B_t)\}}$
 are called the ``discrepancies". We start with two discrepancies at time $t=0$.
 Our challenge is to couple the two permutation processes $A_t$ and $B_t$ so that
 the number of discrepancies is only allowed to decrease (from two to zero).
 Thus we can have at most two discrepancies: one 
 $\begin{pmatrix}
      1    \\
      0  
 \end{pmatrix}$ discrepancy (we will denote it by $d^+_t$)
 and one 
 $\begin{pmatrix}
      0    \\
      1  
 \end{pmatrix}$ discrepancy (we will denote it by $d^-_t$). Here is an example:
 $$\begin{matrix}
    A_t: \quad \dots  & 1 & 0 & 1 & 1 & 0 & 1 & 0 & 0 & \dots   \\
    B_t: \quad \dots  & 1 & 0 & 0 & 1 & 0 & 1 & 1 & 0 & \dots \\
 \phantom{A_t: \quad \dots  }&  &  & \uparrow &  &  &  & \uparrow &  & \\
 \phantom{A_t: \quad \dots  }&  &  & d^+_t &  &  &  & d^-_t &  & 
 \end{matrix}$$
 We recall a similar coupling construction implemented in the recurrent case
 for symmetric exclusion processes.There, whenever the two discrepancies
 happened to be inside the range of a transposition with positive rate, applying the transposition
 to either $A_t$ or $B_t$ we were canceling the discrepancies (see \cite{spitzer74}).
 In our situation, the {\bf tricky} part is that when the two discrepancies happen to be inside the range of
 a permutation from $\Sigma$, applying the permutation to either $A_t$ or $B_t$,
 even if canceling the original two discrepancies, might create new discrepancies.
 This is the challenge that we have to overcome in this subsection.

\subsubsection{Coupling of two-point processes.} \label{conditions}

  Here we will consider three two-point processes $I_t$, $J_t$ and $E_t$ in $S$ with
  the same  initial configuration $x=(x_1,x_2)$ such that $x_1 \not= x_2$.
  We will construct two couplings, one of $I_t$ and $J_t$, and one of $E_t$ and $J_t$.
  First we need to define the processes. 
  
  We assume that the permutation rates $\{ q(\sigma) \}_{\sigma \in \Sigma}$ are known.
  We define $I_t=\{I_1(t), I_2(t)\}$ to be the process consisting of two {\bf independent}
  one-point permutation processes $I_1(t)$ and $I_2(t)$ on $S$, that is two independent 
  one-point permutation processes (random walks) projected on the same space. 
  
  Now, we let $J_t=\{J_1(t), J_2(t)\} \subset S$ be the two-point process that depends on $I_t$ in the
  following way. The initial configuration must be the same: $(I_1(0),I_2(0))=(J_1(0),J_2(0))=x$.
  The above one-point processes $I_1(t)$ and $I_2(t)$ live separate lives. For each of the two
  of them, every $\sigma \in \Sigma$ is enacted  with frequency $q(\sigma)$. The total
  frequency will be $2q(\sigma)$. However, the permutations acting on one of the one-point
  processes will not affect the other. When constructing $J_t$, of all the permutations acting
  on $I_1(t)$ and $I_2(t)$ separately, we will apply to  $J_t$ only those of them that 
  actually displace one of the two random walkers $I_1$ or $I_2$ . Hence, at every moment of time,
  we are waiting for the permutations that contain at least one of the two points ($I_1$ and $I_2$),
  assigning the corresponding $q$-rate to those containing exactly one of them in the range,
  and twice the $q$-rate to those containing both in the range.
   
  Observe, that $I_t$ and $J_t$ are naturally coupled until the ``decoupling" time $T_{dec}$
  when a permutation containing both $J_1(T_{dec}-)$ and $J_2(T_{dec}-)$ occurs 
  (``$t-$" signifies time preceding $t$ such that no changes occur in
  $[t-,t)$ time interval). So $J_1(T_{dec}) \not= J_1(T_{dec}-)$ and 
  $J_2(T_{dec}) \not= J_2(T_{dec}-)$.
  Such permutation should happen before $I_1(t)-I_2(t)$ visits zero for the first time. Thus
  \begin{eqnarray} \label{gg}
  P^x \Big \{ \exists T_{dec}<\infty \text{ s.t. } J_1(T_{dec}) \not= J_1(T_{dec}-) \text{  and  } 
  J_2(T_{dec}) \not= J_2(T_{dec}-) \Big \}
  \end{eqnarray}
   $$\geq P^x \Big \{ \exists t<\infty \text{ s.t. } I_1(t) = I_2(t) \Big \},$$
  where $P^x$ is the probability measure when the corresponding two-point process 
 $I_t$ or $J_t$ (and later the permutation process $E_t$) is at $x \in S^2$ 
  outside the diagonal at time $t=0$.
  We recall that $I_1(t)-I_2(t)$ is recurrent. Hence the left hand side probability
  above is equal to one. As it will be seen soon, this is the primary
  reason why conditions (\ref{I}) and (\ref{II}) are necessary for the coupling construction 
  in \ref{coupling} that follows. 
    
  Now, on the time interval from zero until the decoupling time $T_{dec}$ the process $J_t$
  behaves almost as a two-point permutation process. The only difference being the double
  rates applied to the permutations containing together $J_1$ and $J_2$ in the range at the moment.
  Thus, we find it natural to couple $J_t$ with a two-point exclusion process $E_t=\{E_1(t), E_2(t)\}$
  obeying the same fixed $q$-rates. Lets do that, and on the way clarify the whole construction.
  Define sets $\Sigma_1(t) := \{ \sigma \in \Sigma : I_1(t) \in Range(\sigma) \}$ and
  $\Sigma_2(t) := \{ \sigma \in \Sigma : I_2(t) \in Range(\sigma) \}$. Each permutation in each
  of the two sets occurs with the corresponding $q$-rate, where each permutation in 
  $\Sigma_1(t) \bigcap \Sigma_2(t)$ is counted twice as if two different permutations. 
  Think of $\Sigma_1$ and $\Sigma_2$ as two sets of permutations, of which some are
  identical, but we do not know it and assign separate rates anyways. If the first permutation
  to occur is from $\Sigma_1(t)$, it will act on $I_1$ but not $I_2$, and if it is from $\Sigma_2(t)$, it will act on
  $I_2$, but not $I_1$. No matter to which of the two sets it belongs, the same permutation
  will act on both $J_1$ and $J_2$ even if both are in the range (in the later case the processes
  decouple, and $T_{dec}$ is set to be equal to the action time of such permutation). 
  The same permutation will
  act on both $E_1$ and $E_2$ but only if it comes from $\Sigma_1(t)$ or 
  $\Sigma_2(t) \backslash  \Sigma_1(t) := \{ \sigma \in \Sigma_2(t) : I_2(t) \in Range(\sigma), 
  I_1(t) \not\in Range(\sigma) \}$.
  Of course, $\Sigma_1$ and $\Sigma_2$ evolve after each transformation of $I_t$.
  After decoupling, the processes $I_t$, $J_t$ and $E_t$ evolve independently,
  where $I_t$ is the process consisting of two one-point permutation processes, $E_t$ is a two-point 
  permutation process and $J_t$ is a two-point process where the corresponding $q$-rates 
  are assigned to all
  permutations in $\Sigma$ except for those containing both points $J_1$ and $J_2$ in the range at 
  the moment, assigning the doubled rates to them.

  %We notice that (at least) up until time $T$, $J_t$ evolves according to the original permutation 
  %law except here the rates of permutations that act on both points ($J_1(t)$ and $J_2(t)$) are 
  %being doubled. 
  %Suppose $E_t=\{E_1(t), E_2(t)\}$ is a two-point process that obeys the permutation law, with
  %the same initial configuration \big($(E_1(0),E_2(0))=(I_1(0),I_2(0))=(J_1(0),J_2(0))=x$\big).
  
  For each $\sigma \in \Sigma$, the corresponding Poisson process with frequency
  $q(\sigma)$ can be embedded into a Poisson process with twice the frequency (that is 
  $2q(\sigma)$).  Let $T_{1 \over 2}(\sigma)$ denote the set of jump times for the double-frequency
  Poisson process, then at each point in the time set $T_{1 \over 2}(\sigma)$, the $\sigma$ 
  permutation is either applied to $E_t$ with probability ${1 \over 2}$, or not applied with 
  probability ${1 \over 2}$. When $\sigma \in \Sigma_1 \bigcap \Sigma_2$, that determines 
  whether the permutation comes from $\Sigma_1$ or from $\Sigma_2$.  
  Now, before $E_t$ and $J_t$ decouple, if
  $\sigma \in \Sigma$ and $t \in T_{1 \over 2}(\sigma)$ are such that $E_1(t),E_2(t) \in Range(\sigma)$,
  then $J(t)=\sigma(J(t-))$. 
  %So, at such $t$, the processes decouple with probability ${1 \over 2}$, or not. 
  Thus
  %, since $J(t-)=E(t-) \in Range(\sigma)$,
  \begin{eqnarray} \label{gg1}
  P^x \Big \{ \exists t<\infty \text{ s.t. } t \in T_{1 \over 2}(\sigma) \text{ and }
  E_1(t),E_2(t) \in Range(\sigma) \text{ for some } \sigma \in \Sigma \Big \}  
  \end{eqnarray}
  $$ = P^x \Big \{ \exists t<\infty \text{ s.t. } J_1(t) \not= J_1(t-) \text{  and  } J_2(t) \not= J_2(t-) \Big \}.$$
  At such $t$, either $E(t)=E(t-)$ with probability ${1 \over 2}$ or $E(t)=\sigma \big( E(t-) \big)$.
  In the first case the processes decouple. Since the right hand side of (\ref{gg1}) is equal to
  one in the recurrent case (see (\ref{gg}),   
  $$P^x \Big \{ \exists t<\infty \text{ s.t. } t \in T_{1 \over 2}(\sigma) \text{ and }
  E_1(t),E_2(t) \in Range(\sigma) \text{ for some } \sigma \in \Sigma \Big \} =1$$
  no matter what the starting point $x=(x_1,x_2)$ (s.t. $x_1 \not= x_2$) is.
  So, such $t$ should arrive infinitely often.
  Hence, in the recurrent case, 
  \begin{eqnarray} \label{rec}
  P^x \Big \{ \exists t<\infty \text{ s.t. } E_1(t) \not= E_1(t-) \text{  and  } E_2(t) \not= E_2(t-) \Big \}=1.
  \end{eqnarray}
  It is natural to compare processes $I_t$, $E_t$ and $J_t$ since all three of them coincide 
  up until a certain decoupling time $T_{dec}$.

\subsubsection{The coupling.} \label{coupling}

 We will now try to reconstruct the Spitzer's coupling proof (see \cite{spitzer74}) 
 in the case when conditions (\ref{I}) and (\ref{II}) are satisfied by the permutation 
 process. Lets denote by
 $\Sigma_{cyclic}$ the set of all cyclic permutations in $\Sigma$. We will say that
 a subset $R \subset S$ is a ``range set" if there is a $\sigma \in \Sigma$ with
 $Range(\sigma) = R$.
 Consider a range set $R$.
 Let $$m(R) =\min_{\sigma \in \Sigma : Range(\sigma) = R} \{q(\sigma)\} $$ 
 and $$Z(R)=\sum_{\sigma \in \Sigma : Range(\sigma)=R} q(\sigma).$$
 First, observe that for all range sets $R$ that contain both discrepancies 
 ${\{(A_t \cup B_t) \backslash (A_t \cap B_t)\}}$ at the same time, the sum
 $$z_{d}(t) := \sum_{\footnotesize \begin{matrix} \mbox{range sets }R:\\
                       d^-_t, d^+_t \in R
                       \end{matrix} } Z(R) \leq M_{PL}.$$
 We let the coupled process 
 $\begin{pmatrix}
      A_{t} \\
      B_{t}
  \end{pmatrix}$ evolve according to the following transition rates.
 For each range set $R$ containing both discrepancies at time $t$ we pick a {\bf cyclic} permutation 
 $\sigma_R \in  \Sigma_{cyclic}$ of range $R$ 
 such that $\sigma_R (A_{t})=B_{t}$ (there must be at least one such cyclic permutation). For each range
 set we can order all cyclic permutations, and pick the first one that satisfies the description. Then 
 $$\begin{pmatrix}
      A_{t} \\
      B_{t}
  \end{pmatrix}
  \text{ transforms into }
  \begin{cases}
   &  \begin{pmatrix}
      \sigma_R^2({ A_{t}}) \\
      \sigma_R({B_{t}})
  \end{pmatrix}
  = \begin{pmatrix}
      \sigma_R({B_{t}}) \\
      \sigma_R({B_{t}})
  \end{pmatrix}  
  \text{  with rate } m(R), \\
   \\
   &  \begin{pmatrix}
      \sigma_R^3({ A_{t}}) \\
      \sigma_R^2({B_{t}})
  \end{pmatrix}
  = \begin{pmatrix}
      \sigma_R^2({B_{t}}) \\
      \sigma_R^2({B_{t}})
  \end{pmatrix}  
  \text{  with rate } m(R), \\
  
  & ... \\
  
   &  \begin{pmatrix}
      \sigma_R^{|R|-1}({ A_{t}}) \\
      \sigma_R^{|R|-2}({B_{t}})
  \end{pmatrix}
  = \begin{pmatrix}
      \sigma_R^{-2}({B_{t}}) \\
      \sigma_R^{-2}({B_{t}})
  \end{pmatrix} 
  \text{  with rate } m(R), \\
    \\
    &  \begin{pmatrix}
      \sigma_R ({ A_{t}}) \\
      \sigma_R^{|R|-1}({B_{t}})
  \end{pmatrix}
  =  \begin{pmatrix}
      B_{t} \\
      A_{t}
  \end{pmatrix}   
  \text{  with rate } m(R), \\
    \\
    &  \begin{pmatrix}
      \sigma_R ({ A_{t}}) \\
      \sigma_R ({B_{t}})
  \end{pmatrix}
  =  \begin{pmatrix}
      B_{t} \\
      A_{t}
  \end{pmatrix}   
  \text{  with rate } q(\sigma_R)-m(R), \\
    \\
    &  \begin{pmatrix}
      \sigma_R^2 ({ A_{t}}) \\
      \sigma_R^2 ({B_{t}})
  \end{pmatrix}  
  \text{  with rate } q(\sigma_R^2)-m(R), \\
  
  & ... \\
  
   &  \begin{pmatrix}
      \sigma_R^{|R|-1}({ A_{t}}) \\
      \sigma_R^{|R|-1}({B_{t}})
  \end{pmatrix}
  \text{  with rate } q(\sigma_R^r)-m(R), \\
    \\
    &  \begin{pmatrix}
      \sigma ({ A_{t}}) \\
      \sigma ({B_{t}})
  \end{pmatrix}   
  \text{  with rate } q(\sigma) \text{ if } Range(\sigma)=R \text{ and }\sigma \not= \sigma_R^i, \text{ all } i. \\
 \end{cases}$$
% where $r=r(\sigma_R)$ is 
% the smallest positive integer number such that $\sigma_R^{r+1}= id$. 
 The coupled process 
 $\begin{pmatrix}
      A_{t}\\
      B_{t}
  \end{pmatrix}$ will transform into $\begin{pmatrix}
      \sigma ({ A_{t}}) \\
      \sigma ({B_{t}})
  \end{pmatrix}$ with rate $q(\sigma)$ if $Range(\sigma)$ does not contain both discrepancies.
 %and
  %$$\bar{q}(\sigma_R)= q(\sigma_R)
  %-\sum_{ \sigma \in \Sigma: \sigma^i=\sigma_R \text{ some } i, \text{ }\sigma(A_{t})=B_{t} }
  %{m(Range(\sigma)) \over M_I!} \cdot {Z(Range(\sigma)) \over z_{d}(t)}.$$
 We observe that the rates are well defined.
 %Here, for each such $\sigma_R$,  $\sigma^i=\sigma_R \text{ for some } i$ implies
 %$Range(\sigma_R) \subseteq Range(\sigma)$, and therefore
 %$m(Range(\sigma)) \leq q(\sigma_R)$ (by definition of $m(R)$).
 %Thus\\
 %$\sum_{ \sigma \in \Sigma: \sigma^i=\sigma_R \text{ some } i, \text{ }\sigma(A_{t})=B_{t} }
 %{m(Range(\sigma)) \over M_I!} \cdot {Z(Range(\sigma)) \over z_{d}(t)}$
  %$$\leq q(\sigma_R) \sum_{ \sigma \in \Sigma: \sigma^i=\sigma_R \text{ some } i, \text{ }\sigma(A_{t})=B_{t} }
  %{1 \over M_I!} \cdot {Z(Range(\sigma)) \over z_{d}(t)} $$
  %$$\leq q(\sigma_R) \sum_{U: U \ni \{(A_t \cup B_t) \backslash (A_t \cap B_t)\}}
  %\left({Z(U) \over z_{d}(t)} \cdot  \sum_{\sigma \in \Sigma: Range(\sigma)=U}{1 \over M_I!} \right) 
  %\leq q(\sigma_R)$$ 
  %as $\sum_{U: U \ni \{(A_t \cup B_t) \backslash (A_t \cap B_t)\}}  {Z(U) \over z_{d}(t)} = 1$
  %and there are less then $M_{I}!$ permutations in $\Sigma$ of any given range.
  %The above implies $\bar{q}(\sigma_R) \geq 0$. 
  We also observe that the transformations that
  we have allowed to have non-zero rates do not increase the number of discrepancies.
  Moreover there could be a positive probability of the discrepancies disappearing, in which
  case we let $A_t$ and $B_t$ evolve simultaneously as just a permutation process. The rates
  sum up enabling us to conclude that the above process is a well-defined coupling of 
  processes $A_t$ and $B_t$.

\subsubsection{The coupling is successful. Example.}
  The coupling is successful because, according to (\ref{rec}), if waiting with rate $z_d(t)$ for a permutation
  that contains both discrepancies in its range, though $z_d(t)$ changes with time, we are guaranteed
  a finite holding time. Now, (\ref{I}) and (\ref{II}) imply ${m(R)M_{II}\mathcal{P}(M_{I}) \geq Z(R)}$,
  where $\mathcal{P}(n) = \sum_{k=0}^n \begin{pmatrix}
      n  \\
      k  
  \end{pmatrix} (-1)^k (n-k)! <n!$ denotes the number of permutations of $n>1$ distinct elements
  with all elements displaced (that is
  element $k$ is not in the $k$-th position for all $k \in \{1,2,\dots,n\}$).
  \footnote{ A simple but beautiful Euler's argument shows that 
  $\mathcal{P}(n)=(n-1)(\mathcal{P}(n-1)+\mathcal{P}(n-2))$
  used in one of the many possible derivations of $\mathcal{P}(n)$. Notice that $\mathcal{P}(n)$
  is obviously increasing with $n$.
  Finding the expression for $\mathcal{P}(n)$ is a case of a famous problem, known
  in the history of mathematics by its French name, ``probl\`eme des rencontres". We refer
  the reader to Chapters 3 and 8 of \cite{riordan} for more on the subject.}  
  At the holding time, the discrepancies will cancel with probability
  $$\geq \sum_{\footnotesize \begin{matrix} \mbox{range sets }R:\\
                       d^-_t, d^+_t \in R
                       \end{matrix} } {m(R) \over z_d(t)}
  \geq \sum_{\footnotesize \begin{matrix} \mbox{range sets }R:\\
                       d^-_t, d^+_t \in R
                       \end{matrix} } {Z(R) \over \mathcal{P}(M_{I}) M_{II} z_d(t)}={1 \over \mathcal{P}(M_{I}) M_{II}}.$$ 
  The coupled process will keep arriving to such holding times up until the discrepancies cancel.
  \\
  \\
  {\bf Example.} The author wishes to thank the referee for suggesting the following simple example
  that illustrates how the above coupling works. Let $S=\mathbb{Z}$,\\ 
  $\Sigma=\bigcup_{x \in \mathbb{Z}} 
    \big\{  \sigma_x :=  \overline{(x,x+1,x+2)}, \sigma_x^2=\sigma_x^{-1} \big\}$
  and $q(\sigma_x)=q(\sigma_x^{-1})=q$ for all $x \in \mathbb{Z}$, where $q>0$ is fixed.
  Then one gets $M_{PL}=6q$, $M_{I}=3$ and $M_{II}=1$ (see (\ref{I}) and (\ref{II})).
  Relevant range sets are $R_x=\{x,x+1,x+2\}$ for $x \in \mathbb{Z}$. There $m(R_x)=q$
  and $Z(R_x)=2q$ since $R_x=Range(\sigma_x)=Range(\sigma_x^{-1})$. Suppose the
  discrepancies are for instance at $y$ and $y+1$, say $d^+_t=y$ and $d^-_t=y+1$,
  and the rest of the points around $y$ are occupied in the following way:
 $$\begin{matrix}
    A_t: \quad \dots  & 1 & 1 & 1 & 0 & 0 & 1 & 0 & 0 & \dots   \\
    B_t: \quad \dots  & 1 & 1 & 0 & 1 & 0 & 1 & 0 & 0 & \dots \\
 \phantom{A_t: \quad \dots  }&  & \uparrow & \uparrow & \uparrow &  &  &  &  & \\
 \phantom{A_t: \quad \dots  }&  & y-1 & y & y+1 &  &  \phantom{y+1} &  &  \phantom{y+1} & 
 \end{matrix}$$   
 There are exactly two range sets that contain both discrepancies $d^+_t=y$ and $d^-_t=y+1$, 
 those are $R_{y-1}$ and $R_{y}$. For the range set $R=R_{y-1}$ there is a unique choice of
 $\sigma_R$: $\sigma_R=\sigma_{y-1}^{-1}$ is the only cyclic permutation in $\Sigma$ with 
 range $R$ such that
 $\sigma_R(A_t)=B_t$. Similarly for $R=R_y$, the choice $\sigma_R=\sigma_y$ for $\sigma_R$ is unique. 
 Thus the coupling in \ref{coupling} reads
 $$\begin{pmatrix}
      A_{t} \\
      B_{t}
  \end{pmatrix}
  \text{ transforms into }
  \begin{cases}
   &  \begin{pmatrix}
      \sigma_{y-1}({ A_{t}}) \\
      \sigma_{y-1}^{-1}({B_{t}})
  \end{pmatrix}
  =\begin{matrix}
     & & & & & \\
     & & & & & \\
    \dots & 0 & 1 & 1 & 0 & \dots \\
    \dots & 0 & 1 & 1 & 0 & \dots \\
   \phantom{\dots} & \uparrow & \uparrow & \uparrow &  & \\
   \phantom{\dots} & y-1 & y & y+1 &  &   
 \end{matrix}  
  \text{  with rate } q, \\
   &  \begin{pmatrix}
      \sigma_{y-1}^{-1}({ A_{t}}) \\
      \sigma_{y-1}({B_{t}})
  \end{pmatrix}
  = \begin{matrix}
     & & & & & \\
     & & & & & \\
    \dots & 1 & 0 & 1 & 0 & \dots \\
    \dots & 1 & 1 & 0 & 0 & \dots \\
   \phantom{\dots} & \uparrow & \uparrow & \uparrow &  & \\
   \phantom{\dots} & y-1 & y & y+1 &  &   
 \end{matrix}  
  \text{  with rate } q, \\ 
   &  \begin{pmatrix}
      \sigma_y^{-1}({ A_{t}}) \\
      \sigma_y({B_{t}})
  \end{pmatrix}
  = \begin{matrix}
     & & & & & \\
     & & & & & \\
    \dots & 1 & 0 & 0 & 1 & \dots \\
    \dots & 1 & 0 & 0 & 1 & \dots \\
   \phantom{\dots} & \uparrow & \uparrow & \uparrow &  & \\
   \phantom{\dots} & y-1 & y & y+1 &  &   
 \end{matrix}  
  \text{  with rate } q, \\
   &  \begin{pmatrix}
      \sigma_y({ A_{t}}) \\
      \sigma_y^{-1}({B_{t}})
  \end{pmatrix}
  = \begin{matrix}
     & & & & & \\
     & & & & & \\
    \dots & 1 & 0 & 1 & 0 & \dots \\
    \dots & 1 & 1 & 0 & 0 & \dots \\
   \phantom{\dots} & \uparrow & \uparrow & \uparrow &  & \\
   \phantom{\dots} & y-1 & y & y+1 &  &   
 \end{matrix}  
  \text{  with rate } q.
 \end{cases}$$
  The four permutations that contain both $d^+_t=y$ and $d^-_t=y+1$ have rate
  $z_d(t)=Z(R_y)+Z(R_{y-1})=4q$. If, after waiting with rate $z_d(t)$, the holding time arrives
  (before any changes within $R_y \bigcup R_{y-1}$ occur), the discrepancies will cancel
  with probability equal to ${2q \over z_d(t)}={1 \over 2}$. In general, in all such cases
  when the discrepancies are within distance $\leq 2$ from each other and the holding
  time for all the permutations containing the two discrepancies rings, the probability of
  cancelation of discrepancies should be no less than ${1 \over \mathcal{P}(M_I) M_{II}}={1 \over 2}$
  as $\mathcal{P}(3)=2$. In this example, it will always be equal to ${1 \over 2}$.
    
  The case is obviously recurrent as the difference of corresponding one-point processes 
  $I_1(t)-I_2(t)$ is a recurrent random walk on $\mathbb{Z}$. 
  One can show (see the argument in \ref{conditions}) 
  that for the two-point permutation process $E_t=\{E_1(t), E_2(t)\}$ with the rates given in the
  beginning of the example, the above recurrence implies
  that $E_1(t)$ will come within distance $\leq 2$ of $E_2(t)$ infinitely often insuring that 
  the coupling is successful. 
  
 %\subsection{A different kind of coupling.}
% The examples given in this subsection are interesting because the coupling technique
 %that we apply  here is different from standard (Doeblin) coupling.
  %\textbf{Transpositions.} 
 % \textbf{Triplets (3-cycles).}
  %Consider the case when all the permutations in $\Sigma$ are
  %three-cycles with positive rates. 
  
\subsection{Transient, translation invariant case.}

 We now define the probabilities some of which we already used in the preceeding
 subsections. We let
  $$\bar{g}_2(x) 
  := P^x \Big \{ \exists t<0 \text{ s.t. } E_1(t) \not= E_1(t-) \text{ and } E_2(t) \not= E_2(t-) \Big \},$$
  $$g_2(x) 
  := P^x \Big \{ \exists t<0 \text{ s.t. } I_1(t) = I_2(t) \Big \}$$  
  and
  $$\bar{\bar{g}}_2(x) 
  := P^x \Big \{ \exists t<0 \text{ s.t. } J_1(t) \not= J_1(t-) \text{ and } J_2(t) \not= J_2(t-) \Big \},$$
  where $P^x$ is again the probability measure when the corresponding two-point process 
  $E_t$, $I_t$ or $J_t$ is at $x \in S^2$ at time $t=0$.
  Therefore (\ref{gg}) is equivalent to $$\bar{\bar{g}}_2(x) \geq g_2(x).$$ 
  Moreover, by construction, $\bar{\bar{g}}_2 (x) \geq \bar{g}_2(x) \geq g_2(x)$. The equality 
  (\ref{gg1}) implies $\bar{g}_2(x) \geq {1 \over 2} \bar{\bar{g}}_2(x)$,
  and one similarly obtains ${1 \over M_{II} \mathcal{P}(M_I)} \bar{\bar{g}}_2(x) \leq g_2(x)$,
  where as before, $\mathcal{P}(N)$ denotes the number of permutations of $N$ elements 
  with all elements displaced. 
  Hence, taking all the above inequalities together, we conclude that 
  \begin{eqnarray} \label{similar}
  g_2(x) \sim \bar{g}_2(x) \sim \bar{\bar{g}}_2(x).
  \end{eqnarray}

  Now, let
  $$T_n =\{ x=(x_1,\dots,x_n)\in S^n \text{ :  } x_i \not= x_j \text{  for all  }i \not= j \}.$$  
  Let $S_t$, $U_t$  and $V_t$ be the semigroups of respectively $E_t$, $I_t$ and $J_t$.
  If we let $E_t$ be the $n$ points permutation process and generalize $I_t$ to be the
  corresponding $n$ points process, where each point moves independently of the others
  as a one-point permutation process, then we can redefine
  $$\bar{g}_n(x) 
  := P^x \Big \{ \exists t<0 \text{ s.t. } E_i(t) \not= E_i(t-) \text{ and }
  E_j(t) \not= E_j(t-) \text{ for some } i \not= j \in \{1,\dots,n\} \Big \}$$
  and
  $$g_n(x) := P^x \Big \{ \exists t<0 \text{ s.t. } I_t=(I_1(t),\dots,I_2(t)) \not\in T_n \Big \}.$$ 
  The properties of $g_n$ were thoroughly studied before (see for example \cite{liggett1}).
  In particular, for $x=(x_1,\dots,x_n) \in S^n$,
  \begin{eqnarray*} %\label{g2g}
  g_n(x) \leq \sum_{1 \leq i<j \leq n} g_2(x_i, x_j) \leq 
  \begin{pmatrix}
      n\\
      2  
\end{pmatrix} g_n(x).
  \end{eqnarray*}
  Thus, redoing the above estimates for a general $n$, one gets 
  \begin{eqnarray} \label{nsimilar}
  \bar{g}_n(x) \sim g_n(x) \sim \sum_{1 \leq i<j \leq n} g_2(x_i, x_j)
  \sim \sum_{1 \leq i<j \leq n} \bar{g}_2(x_i, x_j)
  \end{eqnarray}
  
\subsubsection{Case $n=2$}
  By following Liggett's proof (see Theorem 1.44 in Chapter  VIII of
   \cite{liggett1}, \cite{liggett3} and \cite{liggett3a}) for transient symmetric exclusion process, we observe that,
  by construction, if $f$ is a function such that $0\leq f \leq 1$, then
  \begin{eqnarray*}
  |V_t f(x)-U_t f(x)| \leq \bar{\bar{g}}_2(x), \text{\qquad} x \in T_2.
  \end{eqnarray*}  
  Here $J_t$ and $I_t$ agree until the first time $t$ such that $J_1(t) \not=J_1(t-)$
  and $J_2(t) \not=J_2(t-)$. Now, $J_t$ agrees with $E_t$ up until at least such $t$.
  Thus
  \begin{eqnarray*}
  |V_t f(x)-S_t f(x)| \leq \bar{\bar{g}}_2(x), \text{\qquad} x \in T_2.
  \end{eqnarray*}
  and, by (\ref{similar}),
  \begin{eqnarray} \label{ineqT1}
  |S_t f(x)-U_t f(x)| \leq 2\bar{\bar{g}}_2(x) \leq  \bar{g}_2(x), \text{\qquad} x \in T_2.
  \end{eqnarray}
  
  Suppose $f$ is also symmetric on $T_2$, and $S_t f=f$ for all $t \geq 0$. It can be
  extended to all of $S^2$ by setting $f=0$ on $T_2^c :=S^2 \backslash T_2$. Then, by (\ref{ineqT1}),
  \begin{eqnarray} \label{ineqT2} 
  |f(x)-U_t f(x)| \leq \bar{g}_2(x), \text{\qquad} x \in S^2
  \end{eqnarray}
  as $\bar{g}_2 := 1$ on  $T_2^c$.
  
  We refer the reader to \cite{liggett1} for the proof of 
  $$\lim_{t\rightarrow \infty} U_t g_n(x)=0, \text{\qquad} x\in S^n.$$
  Thus, by (\ref{similar}),
  \begin{eqnarray} \label{limit} 
  \lim_{t\rightarrow \infty} U_t \bar{g}_2(x)=0, \text{\qquad} x\in S^2.
  \end{eqnarray}
  The inequality (\ref{ineqT2}) implies
  $$ |U_s f(x)-U_{s+t} f(x)| \leq U_s\bar{g}_2(x), \text{\qquad} x \in S^2,$$
  where, by (\ref{limit}), the right hand side goes to zero. So, the limit of $U_s f$ exists
  and is $U_t$-harmonic, whence it is a constant
  $$\lim_{t\rightarrow \infty} U_t f(x)=C, \text{\qquad} x\in S^2.$$
  Thus (\ref{ineqT2}) implies
  \begin{eqnarray*}
  |f(x)-C| \leq \bar{g}_2(x), \text{\qquad} x \in S^2.
  \end{eqnarray*}
  Since we know that $S_t f=f$,
  \begin{eqnarray} \label{Sff}
  |f(x)-C| =|S_t f(x) -C| \leq S_t \bar{g}_2(x), \text{\qquad} x \in T_2.
  \end{eqnarray}

    \textbf{Three-cycles.} If we allow only transpositions and three-cycles then the situation will be
   much simpler. First consider the case when $\Sigma$ contains only three-cycles. So, we only have to
   consider the permutations $\sigma_z$, indexed by $z \not= x_1$  or $x_2$ in $S$ such
   that $\sigma_z: \mbox{  } z \rightarrow x_1 \rightarrow x_2 \rightarrow z$, as well as 
   $\sigma_z^{-1}$. Let $\Omega$, $\mathsf{U}$ and $\mathsf{V}$ be the generators of 
   the corresponding semigroups  $S_t$, $U_t$ and $V_t$.
   For a cylinder function $h: S \times S \rightarrow \mathbb{R}$ and $x=(x_1, x_2) \in S^2$,
  $$(\mathsf{U}-\mathsf{V}) h(x) = \sum_{\sigma: x_1, x_2 \in Range(\sigma)}
      q(\sigma) \Big[ h(\sigma(x_1), x_2)+h(x_1, \sigma(x_2)) -2h(\sigma(x_1), \sigma(x_2)) \Big]$$
  and
  $$(\mathsf{V}-\Omega) h(x) = \sum_{\sigma: x_1, x_2 \in Range(\sigma)}
      q(\sigma) \Big[h(\sigma(x_1), \sigma(x_2)) -h(x_1, x_2) \Big].$$
  Thus
  \begin{eqnarray} \label{U-S}
  (\mathsf{U}-\Omega) h(x) = \sum_{\sigma: x_1, x_2 \in Range(\sigma)}
      q(\sigma) \Big[ h(\sigma(x_1), x_2)+h(x_1, \sigma(x_2)) 
      -h(\sigma(x_1), \sigma(x_2)) -h(x_1, x_2) \Big].
   \end{eqnarray} 
      
Here taking the portion of the sum in (\ref{U-S}) corresponding to the
   three-cycles $\sigma_z$ and $\sigma_z^{-1}$ we obtain the following equality:
   $$(\mathsf{U}-\Omega) h(x) = \sum_{\sigma: Range(\sigma)=\{z, x_1, x_2 \} }
      q(\sigma) \Big[ h(\sigma(x_1), x_2)+h(x_1, \sigma(x_2)) 
      -h(\sigma(x_1), \sigma(x_2)) -h(x_1, x_2) \Big]$$
   $$=q(\sigma_z) \Big[ h(x_2, x_2)+h(x_1, z)  -h(x_2, z) -h(x_1, x_2) \Big]
      + q(\sigma_z) \Big[ h(z, x_2)+h(x_1, x_1) -h(z, x_1) -h(x_1, x_2) \Big]$$
   $$=q(\sigma_z) \Big[ h(x_1, x_1)+h(x_2, x_2)-2h(x_1, x_2) \Big].$$
   A bounded symmetric function $F$ on $S^2$ is said to be positive definite 
   if
   \begin{eqnarray} \label{positive}
   \sum_{u_1,u_2 \in S} F(u_1,u_2) \beta(u_1) \beta(u_2) \geq 0
   \end{eqnarray}
   whenever $\sum_{u \in S}  |\beta(u)| < \infty$ and $\sum_{u \in S} \beta(u)=0$.
   A bounded symmetric function $F$ on $S^n$ is said to be positive definite 
   if it is a positive definite function of each pair of its variables.     
   Now, $h(x)=U_s g_2(x)$ is positive definite (see the proof of Lemma 1.23 in Chapter  VIII of
   \cite{liggett1}). Taking $\beta(u)=\begin{cases}
      +1 \text{ if } u=x_1 \\
      -1 \text{ if } u=x_2 \\
      0 \text{ otherwise }
   \end{cases}$
   in (\ref{positive}) we conclude that $(\mathsf{U}-\Omega)U_s g_2(x) \geq 0$. Thus
   \begin{eqnarray}  \label{UgSg}
   S_t g_2(x) \leq U_t g_2(x).
   \end{eqnarray}
   follows from the integration by parts formula for semigroups
   $$U_t-S_t=\int_0^t S_{t-s}(\mathsf{U}-\Omega)U_s ds.$$   
   (\ref{UgSg}) together with (\ref{Sff}) and (\ref{similar}) complete the argument in case when we only allow 
   three-cycles.  The proof can be easily extended to allow $\Sigma$ to include
   both transpositions and three-cycles, by incorporating the proof of
   Proposition 1.7 in Chapter  VIII of \cite{liggett1}. 
   
   \textbf{For the general case} the inequalities like (\ref{UgSg}) are hard to prove. However
    (\ref{UgSg}) is stronger than what we really need.
    
    By transience, $\lim_{x \rightarrow  \infty} g_2(0,x)=0$, $x \in S$. This together with (\ref{similar}) imply \\
    ${\lim_{x \rightarrow  \infty}\bar{g}_2(0,x)=0}$. So, for any $\epsilon >0$ $\exists R_{\epsilon}$
    such that $\bar{g}_2(0,x) \leq \epsilon$ whenever $|x|>R_{\epsilon}$. Now , we claim that there
    exists $\Delta <1$ such that
    $\bar{g}_2(x_1,x_2)=\bar{g}_2(0,x_2-x_1) \leq \Delta$ for all $x_1, x_2 \in S$. To prove this, we
    consider any $\epsilon \in (0,1)$, say $\epsilon={1 \over 2}$, and denote $R=R_{1 \over 2}$.
    We only need to prove that $\bar{g}_2(0,x)<1$ whenever $|x| \leq R$, $x \in S$. Suppose there is a point 
    $x$ inside the ball $B_R$ of radius $R$ around the origin such that $\bar{g}_2(0,x)=1$. If there is
    a permutation $\sigma_1$ of positive rate with $\sigma_1(x) \in B_R^c$ and $0 \not\in Range(\sigma_1)$,
    then
    $$1-\bar{g}_2(0,x) \geq (1-\bar{g}_2(0,\sigma_1(x)) P^{(0,x)}(\eta_t=(0,\sigma_1(x))$$
    for small $t$ such that
    $$P^{(0,x)}(\eta_t=(0,\sigma_1(x)) \geq t q(\sigma_1) e^{-4 M_{PL}T} >0,$$
    where the RHS signifies the case when $\sigma_1$ is the only permutation containing $0$, $x$
    and/or $\sigma_1(x)$ in its range that acts within the interval $[0,t]$ (we recall that $M_{PL}$ comes
   from the permutation law settings, see (\ref{PL}) ).
    
    Thus $\exists \Delta_1 <1$ such that $\bar{g}_2(0,x) \leq \Delta_1$ whenever 
    $$x \in B_R^c \cup \{ x \in B_R \text{  :  } \exists \sigma_1 \in \Sigma \text{  s.t. } 
           \sigma_1(0)=0, \text{  } \sigma_1(x) \in B_R^c \}.$$
           
     Similarly, since there are finitely many points of $S$ inside $B_R$,
      $\exists \Delta <1$ such that ${\bar{g}_2(0,x) \leq \Delta}$ whenever 
     $$x \in B_R^c \cup \{ x \in B_R \text{  :  }\exists k \geq 1, \sigma_1,...,\sigma_k \in \Sigma \text{  s.t. } 
           \sigma_1(0)=...=\sigma_k(0)=0, \text{  } 
           \sigma_k \circ \sigma_{k-1}\circ ... \circ \sigma_1(x) \in B_R^c \}.$$
    By irreducibility assumption, the above set is all of $S$, proving the claim.
    Thus $\forall M>0$,  $P^{(0,x)}(|E_1(t)-E_2(t)|\leq M \text{  i.o.})=0$ and
    $|E_1(t)-E_2(t)| \rightarrow +\infty$ as $t \rightarrow \infty$.
    Thus $\lim_{x \rightarrow  \infty} \bar{g}_2(0,x)=0$ implies
    $$\lim_{x \rightarrow  \infty} S_t \bar{g}_2(0,x)=0.$$
    Hence, by (\ref{Sff}), if $S_t f =f$ then $f(x)$ is a constant for all $x \in T_2$,
    e.g. a bounded harmonic function for the transient permutation process $E_t$ is constant 
    for all sets of cardinality $n=2$, proving 
    Theorem \ref{harmonic} in this case.
    
\subsubsection{General $n$}
  The proof that, if $f$ is a bounded symmetric function on $T_n$, and if $S_t f =f$, then
  \begin{eqnarray} \label{nSff}
  |f(x)-C| =|S_t f(x) -C| \leq S_t \bar{g}_n(x), \text{\qquad} x \in T_n
  \end{eqnarray}  
  for some constant $C$ is the same for general $n$ as in case when $n=2$.
  However, here we do not have to do the rest of the computations again.
  Since $\lim_{x \rightarrow  \infty} S_t \bar{g}_2(0,x)=0$ for all $x \not= 0$ in $S$,
  (\ref{nsimilar}) implies that the right side of (\ref{nSff}) goes to zero. Thus, for all integer $n \geq 2$,
  a bounded harmonic function for the transient permutation process $E_t$ must be constant 
  for all sets of cardinality $n$. Theorem \ref{harmonic} is proved.

\section{General case: shift invariant stationary\\ measures}
 Once again, we assume that the conditions (\ref{I}) and (\ref{II}) are satisfied,
 though, as it was mentioned in the previous section, it is possible to obtain
 some of the same results with slightly weaker conditions.
 
 Let  $\mathcal{S}$ again denote the class of the shift invariant probability measures on $\{0,1\}^S$.
 In this section we will prove the following important
 \begin{thm} \label{ISe}
 For the general permutation process, 
 $(\mathcal{I} \cap \mathcal{S})_e=\{ \nu_{\rho}: 0 \leq \rho \leq 1\}$.
 \end{thm}

\subsection{Modifying the coupling} 
 First we have to modify the coupling of two permutation processes $A_t$ and $B_t$ on $S$, where
 now we are not constrained to only two discrepancies at time $t=0$. We should find the way
  of coupling the processes so that the number of discrepancies is at least not 
 increasing with time. We will adapt the following (generally accepted) notation:
 for two configurations $\eta$ and $\zeta$ in $\{0,1\}^S$, we say that $\eta \leq \zeta$ if
 $$\eta(x) \leq \zeta(x) \mbox{  for every } x \in S.$$ We will say that $\eta \leq \zeta$
 \textbf{on} a subset $S_{sub}$ of $S$ if $\eta(x) \leq \zeta(x)$ for every $x \in S_{sub}.$
 
 At  a given time $t$, for every range set $R$ in $S$, there must be at least one $\sigma_R \in \Sigma_{cyclic}$ 
 of range $R$ (i.e. $Range(\sigma_R)=R$,) such that either $\sigma_R(A_{t}) \geq B_{t}$ 
 \textbf{on} $R$ or $\sigma_R(B_{t}) \geq A_{t}$ \textbf{on} $R$. 
 In the case when 
 $$|\{x \in R: A_t(x)=1\}| \geq |\{x \in R: B_t(x)=1\}|,$$
 we can only pick $\sigma_R$ so that $\sigma_R(A_{t}) \geq B_{t}$. Then we let 
 $\begin{pmatrix}
      A_{t} \\
      B_{t}
  \end{pmatrix}$ transform into
 either  $\begin{pmatrix}
      \sigma_R(A_t) \\
      B_t  
 \end{pmatrix}$,
  $\begin{pmatrix}
      \sigma_R^2(A_t) \\
      \sigma_R(B_t)  
 \end{pmatrix}$,
  $\begin{pmatrix}
      \sigma_R^3(A_t) \\
      \sigma_R^2(B_t)  
 \end{pmatrix}$, $\dots$, or
  $\begin{pmatrix}
      \sigma_R^{|R|}(A_t) \\
      \sigma_R^{|R|-1}(B_t)  
  \end{pmatrix}=\begin{pmatrix}
      A_t \\
      \sigma_R^{-1}(B_t)  
  \end{pmatrix}$ with rate $m(R)$ each, where $m(R)$ is defined as in \ref{coupling}.
  For all permutations  $\sigma \in \Sigma$ of range $R$, we will apply  $\begin{pmatrix}
      \sigma \\
      \sigma  
 \end{pmatrix}$ with the remaining rates: 
 $\begin{pmatrix}
      A_{t} \\
      B_{t}
  \end{pmatrix}$ transforms into 
  $\begin{pmatrix}
      \sigma(A_{t}) \\
      \sigma(B_{t})
  \end{pmatrix}$ with rate $=q(\sigma)-m(R)$ if $\sigma=\sigma_R^i$ for some $i \in \{1,...,|R|-1\}$,
  and with rate $=q(\sigma)$ if otherwise.   The case when $\sigma(B_{t}) \geq A_{t}$ on $Range(\sigma)$ 
  is dealt with symmetrically. The way we select $\sigma_R$ among the cyclic permutations
  of range $R$ is by initially ordering all cyclic permutations of range $R$, and at every time selecting 
  the one of highest order s.t. $\sigma_R(A_{t}) \geq B_{t}$.  
  It is {\bf important} that the ordering of all the cyclic permutations of range $R$ should be done parallel to
  ordering of cyclic permutations of range $R+y$ for each $y \in S$, i.e. $\sigma_R$ selected for
  $\begin{pmatrix}
      A_t(x) \\
      B_t(x)
  \end{pmatrix}
  =\begin{pmatrix}
      \eta(x) \\
      \zeta(x)
  \end{pmatrix}$ for all $x \in S$ should be the $(-y)$-shift of $\sigma_{R+y}$ selected for 
    $\begin{pmatrix}
      A_t(x) \\
      B_t(x)
  \end{pmatrix}
  =\begin{pmatrix}
      \eta(x-y) \\
      \zeta(x-y)
  \end{pmatrix}$ for all $x \in S$.
  We \textbf{observe} that the number of discrepancies here can only decrease.
  
  We will denote by $\mathcal{I}^*$ the class of stationary distributions for the coupled
  process, and by $\mathcal{S}^*$ we will denote the class of translation invariant distributions
  for the coupled process. We will also write $\mathcal{I}_e^*$ for the set of all the extreme
  points of $\mathcal{I}^*$, and $(\mathcal{I}^* \cap \mathcal{S}^*)_e$ for the set of all the
  extreme points of $(\mathcal{I}^* \cap \mathcal{S}^*)$
  Let $\nu^*$ be the measure on $\{0,1\}^S \times \{0,1\}^S$ with the marginals
  $\nu_1$ and $\nu_2$.  Our next theorem is a case of Theorem 2.15 in Chapter III of 
  \cite{liggett1}.
  \begin{thm} \label{star}
  (a) If $\nu^*$ is in $\mathcal{I}^*$, then its marginals are in $\mathcal{I}$.\\
  (b) If $\nu_1, \nu_2 \in \mathcal{I}$, then there is a $\nu^* \in \mathcal{I}^*$ with marginals
  $\nu_1$ and $\nu_2$.\\
  (c) If $\nu_1, \nu_2 \in \mathcal{I}_e$, then the $\nu^*$ in part (b) can be taken to be in
  $\mathcal{I}_e^*$.\\
  (d) In parts (b) and (c), if $\nu_1 \leq \nu_2$, then $\nu^*$ can be taken to concentrate 
  on $\{\eta \leq \zeta \}$.\\
  (e) In the translation invariant case, parts (a)-(d) hold if $\mathcal{I}$ and $\mathcal{I}^*$
  are replaced by $(\mathcal{I} \cap \mathcal{S})$ and $(\mathcal{I}^* \cap \mathcal{S}^*)$
  respectively.
  \end{thm}
 
\subsection{Case $\nu^* \in (\mathcal{I}^* \cap \mathcal{S}^*)$: 
                     the two types of discrepancies do not coexist}
 For permutations $\sigma_1, \sigma_2 \in S$ of a given range $R$, let 
 $q^*(\sigma_1, \sigma_2; \eta(R), \zeta(R))$ denote the rate of the newly defined coupled process 
 assigned to
 $\begin{pmatrix}
      \sigma_1 \\
      \sigma_2  
 \end{pmatrix}$ transformation if given the values
 $\begin{pmatrix}
      \eta(R) \\
      \zeta(R)  
 \end{pmatrix}=\left \{ \begin{pmatrix}
      \eta(x) \\
      \zeta(x)  
 \end{pmatrix} \mbox{ for all } x \in R \right \}$. 
 We also let $S^*(t)$ denote the semigroup of the coupled process.
 The following definition will be useful as we proceed:
 \begin{Def*}
 For $\sigma \in \Sigma$ and $x \in Range(\sigma)$, the subset
 $$O(\sigma, x)=\{ \sigma^i (x): \quad i=0,1,\dots \}$$
 of $Range(\sigma)$ is called the {\bf orbit} of $x$ under $\sigma$.
 \end{Def*}
 
 \begin{thm} \label{dplus}
 If $\nu^* \in (\mathcal{I}^* \cap \mathcal{S}^*)$, then
 $$\nu^* \big\{(\eta, \zeta)\text{ : } \eta(u)=\zeta(v)=0,\text{ } \zeta(u)=\eta(v)=1 \big\}=0$$
 for every $x$ and $y$ in $S$.
 \end{thm}
 \textit{Proof:}
 Here we reconstruct a clever trick from the theory of exclusion processes.
 If the coupled measure $\nu^* \in (\mathcal{I}^* \cap \mathcal{S}^*)$ then
 \begin{eqnarray} \label{SI}
 0 & = &
 {d \over dt} \nu^* S^*(t)\{(\eta, \zeta)\text{ : }\eta(x) \not= \zeta(x) \} \Big|_{t=0}\\
 \nonumber \\
 & = & \sum_{\footnotesize \begin{matrix} \mbox{range sets }R:\\
                       x \in R
                       \end{matrix} }
          \sum_{\footnotesize \begin{matrix} 
                       \sigma \in \Sigma:\\
                       Range(\sigma)=R
                       \end{matrix} } 
          \sum_{\footnotesize \begin{matrix}
                      \ddot{\eta},\ddot{\zeta} \in \{0,1\}^R:\\
                      \ddot{\eta}(x) = \ddot{\zeta}(x),\\
                      \ddot{\eta}(\sigma^{-1}(x)) \not= \ddot{\zeta}(\sigma^{-1}(x))    
                      \end{matrix} } q^*(\sigma, \sigma; \ddot{\eta}, \ddot{\zeta}) 
 \cdot \nu^* \Big\{ \begin{pmatrix} \eta \\ \zeta \end{pmatrix}: 
        \begin{pmatrix} \eta(R) \\ \zeta(R) \end{pmatrix} = \begin{pmatrix} \ddot{\eta} \\ \ddot{\zeta} \end{pmatrix} 
          \Big\}     \nonumber \\ 
 \nonumber \\ 
 &-&  \sum_{\footnotesize \begin{matrix} \mbox{range sets }R:\\
                       x \in R
                       \end{matrix} }
          \sum_{\footnotesize \begin{matrix} 
                       \sigma \in \Sigma:\\
                       Range(\sigma)=R
                       \end{matrix} } 
          \sum_{\footnotesize \begin{matrix}
                      \ddot{\eta},\ddot{\zeta} \in \{0,1\}^R:\\
                      \ddot{\eta}(x) \not= \ddot{\zeta}(x),\\
                      \ddot{\eta}(\sigma^{-1}(x)) = \ddot{\zeta}(\sigma^{-1}(x))    
                      \end{matrix} } q^*(\sigma, \sigma; \ddot{\eta}, \ddot{\zeta}) 
 \cdot \nu^* \Big\{ \begin{pmatrix} \eta \\ \zeta \end{pmatrix}: 
        \begin{pmatrix} \eta(R) \\ \zeta(R) \end{pmatrix} = \begin{pmatrix} \ddot{\eta} \\ \ddot{\zeta} \end{pmatrix} 
          \Big\}     \nonumber \\ 
 \nonumber \\ 
 &+&  \sum_{\footnotesize \begin{matrix} \mbox{range sets }R:\\
                       x \in R
                       \end{matrix} } 
          \sum_{\footnotesize \begin{matrix}
                      \ddot{\eta},\ddot{\zeta} \in \{0,1\}^R:\\
                      |\{x \in R: \ddot{\eta}(x) =1\}| \\
                      \geq |\{x \in R: \ddot{\zeta}(x) =1\}|    
                      \end{matrix} } 
 {\footnotesize m(R) \big[ D(\sigma_R(\ddot{\eta}), \ddot{\zeta})-D(\ddot{\eta}, \ddot{\zeta}) \big]
 \cdot \nu^* \Big\{ \begin{pmatrix} \eta \\ \zeta \end{pmatrix}: 
        \begin{pmatrix} \eta(R) \\ \zeta(R) \end{pmatrix} = \begin{pmatrix} \ddot{\eta} \\ \ddot{\zeta} \end{pmatrix} 
          \Big\} }     \nonumber \\ 
 \nonumber \\ 
 &+&  \sum_{\footnotesize \begin{matrix} \mbox{range sets }R:\\
                       x \in R
                       \end{matrix} } 
          \sum_{\footnotesize \begin{matrix}
                      \ddot{\eta},\ddot{\zeta} \in \{0,1\}^R:\\
                      |\{x \in R: \ddot{\eta}(x) =1\}| \\
                      \geq |\{x \in R: \ddot{\zeta}(x) =1\}|    
                      \end{matrix} } 
 {\footnotesize m(R) \big[ D(\ddot{\eta}, \sigma_R(\ddot{\zeta}))-D(\ddot{\eta}, \ddot{\zeta}) \big]
 \cdot \nu^* \Big\{ \begin{pmatrix} \eta \\ \zeta \end{pmatrix}: 
        \begin{pmatrix} \eta(R) \\ \zeta(R) \end{pmatrix} = \begin{pmatrix} \ddot{\eta} \\ \ddot{\zeta} \end{pmatrix} 
          \Big\} },       \nonumber 
  \end{eqnarray}
  where for each a range set $R$ and configuration $\ddot{\eta},\ddot{\zeta} \in \{0,1\}^R$
  of the coupled process on $R$, $\sigma_R$ is uniquely defined.  Also $D^+(\ddot{\eta},\ddot{\zeta})$ is the number of 
  $\begin{pmatrix}
      1\\
      0
  \end{pmatrix}$ discrepancies of
  $\begin{pmatrix}
      \ddot{\eta}\\
      \ddot{\zeta} 
 \end{pmatrix}$,
 $D^-(\ddot{\eta},\ddot{\zeta})$ is the number of 
  $\begin{pmatrix}
      0\\
      1
  \end{pmatrix}$ discrepancies of
  $\begin{pmatrix}
      \ddot{\eta}\\
      \ddot{\zeta} 
 \end{pmatrix}$
 and 
 $$D(\ddot{\eta},\ddot{\zeta}) := D^+(\ddot{\eta},\ddot{\zeta})+D^-(\ddot{\eta},\ddot{\zeta})$$
 is the total number of discrepancies on $R$; $\sigma_R(\ddot{\eta})$ above denotes the
 the disposition of the particles that we get after applying permutation $\sigma_R$ to the original
 $\ddot{\eta}$: $\sigma_R(\ddot{\eta})(x) := \ddot{\eta}(\sigma_R^{-1}(x))$ for all $x \in R$, 
 $\sigma_R(\ddot{\zeta})$ is defined by analogy.
 
 Now, here is some explanation. The third sum on the right hand side (RHS) of (\ref{SI}) represents the
 contribution to the derivative by all transformations 
 $\begin{pmatrix}
      \sigma_R^i \\
      \sigma_R^{i-1}  
 \end{pmatrix}$ whenever $$|\{x \in R: \eta(x)=1\}| \geq |\{x \in R: \zeta(x)=1\}|$$ (equivalently
 $\sigma_R(\eta) \geq \zeta$ on $R$). Symmetrically, the fourth sum on the RHS of (\ref{SI}) represents the
 contribution to the derivative by all transformations 
 $\begin{pmatrix}
      \sigma_R^{i-1} \\
      \sigma_R^i  
 \end{pmatrix}$ whenever $|\{x \in R: \zeta(x)=1\}| \geq |\{x \in R: \eta(x)=1\}|$ (equivalently
 $\sigma_R(\zeta) \geq \eta$ on $R$).

 Now, lets show that the third sum is correct. We fix a range set $R$ that contains $x$. Notice that
 since $\sigma_R$ is cyclic, for each $y \in R$ there is a unique corresponding $i \in \{0,1,...,|R|-1 \}$
 such that $y=\sigma_R^{-i}(x)$. So $\sigma_R^{i+1}(\ddot{\eta})(x)=\ddot{\eta}(\sigma_R^{-(i+1)}(x))
 =\ddot{\eta}(\sigma_R^{-1})(y)=\sigma_R(\ddot(\eta))(y)$ and similarly
 $\sigma_R^i(\ddot(\zeta))(x)=\ddot{\zeta}(y)$.
 Then counting all contributions to the derivative in (\ref{SI}) by
  $\begin{pmatrix}
      \sigma_R^{i+1} \\
      \sigma_R^i  
 \end{pmatrix}$ for all values of $i \in \{0,1,...,|R|-1 \}$ given that
 ${ \begin{pmatrix} \eta(R) \\ \zeta(R) \end{pmatrix} = \begin{pmatrix} \ddot{\eta} \\ \ddot{\zeta} \end{pmatrix} }$, 
 one obtains product of 
 $\nu^* \Big\{ \begin{pmatrix} \eta \\ \zeta \end{pmatrix}: 
        \begin{pmatrix} \eta(R) \\ \zeta(R) \end{pmatrix} = \begin{pmatrix} \ddot{\eta} \\ \ddot{\zeta} \end{pmatrix} 
          \Big\}$ with
 $$\sum_{i=0}^{|R|-1} m(R) \Big[ \mathbf{1}_{\{\sigma_R^{i+1}(\eta)(x) \not= \sigma_R^i(\zeta)(x) \}}
     -\mathbf{1}_{\{\eta(x) \not= \zeta(x) \}} \Big]
     =m(R) \sum_{y \in R} \Big[ \mathbf{1}_{\{\sigma_R(\eta)(y) \not= \zeta(y) \}}
     -\mathbf{1}_{\{\eta(x) \not= \zeta(x) \}} \Big]$$
  $$=m(R) \left[ D(\sigma_R(\ddot{\eta}), \ddot{\zeta})-|R|\mathbf{1}_{\{\eta(x) \not= \zeta(x) \}} \right],$$
  where $\mathbf{1}_{\{\lambda \not= \mu\}} :=
  \begin{cases}
      1\text{ }, \\
      \lambda \not= \mu \text{  otherwise}.
  \end{cases}$ 
  
  Next step is to consider all the shifts $R_{x-z}=\{R+x-z\}$ of $R$ for all $z \in R$ together with
  the corresponding shifts $\begin{pmatrix} \ddot{\eta}^{x-z} \\ \ddot{\zeta}^{x-z} \end{pmatrix}$
  of configurations $\ddot{\eta},\ddot{\zeta} \in \{0,1\}^R$. Since $\nu^*$
  is shift invariant, the contribution to the derivative in (\ref{SI}) coming from
  all transitions   $\begin{pmatrix}
      \sigma_{R_{x-z}}^{i+1} \\
      \sigma_{R_{x-z}}^i  
  \end{pmatrix}$ for all values of $i \in \{0,1,...,|R|-1 \}$ and $z \in R$ when each time given that
 ${ \begin{pmatrix} \eta(R_{x-z}) \\ \zeta(R_{x-z}) \end{pmatrix} 
 = \begin{pmatrix} \ddot{\eta}^{x-z} \\ \ddot{\zeta}^{x-z} \end{pmatrix} }$,
 is equal to
  $$\nu^* \Big\{ \begin{pmatrix} \eta \\ \zeta \end{pmatrix}: 
        \begin{pmatrix} \eta(R) \\ \zeta(R) \end{pmatrix} = \begin{pmatrix} \ddot{\eta} \\ \ddot{\zeta} \end{pmatrix} 
          \Big\} \cdot m(R) \sum_{z \in R}  
          \left[ D(\sigma_R(\ddot{\eta}), \ddot{\zeta})-|R|\mathbf{1}_{\{\eta(z) \not= \zeta(z) \}} \right]$$
  $$=|R| \cdot m(R) \big[ D(\ddot{\eta}, \sigma_R(\ddot{\zeta}))-D(\ddot{\eta}, \ddot{\zeta}) \big]
 \cdot \nu^* \Big\{ \begin{pmatrix} \eta \\ \zeta \end{pmatrix}: 
        \begin{pmatrix} \eta(R) \\ \zeta(R) \end{pmatrix} = \begin{pmatrix} \ddot{\eta} \\ \ddot{\zeta} \end{pmatrix} 
          \Big\}.$$
   Now, the above is the total contribution corresponding to $|R|$ shifts of $R$ that still contain $x$ and
   respective shifts of $\begin{pmatrix} \ddot{\eta} \\ \ddot{\zeta} \end{pmatrix}$. 
   Hence we can count in
   ${1 \over |R|}$-th fraction of the total each time, thus verifying the correctness of the third sum on
   the RHS of (\ref{SI}). 
  
 %So, for each $y \in R$, the pair
  %$\begin{pmatrix}
     % \sigma_R(\eta(y))\\
      %\zeta(y) 
 %\end{pmatrix}$
 %visits site $x$ exactly once
 %(e.g. there is exactly one $i$, $0 \leq i \leq |R|-1$, such that
   %$\begin{pmatrix}
      %\sigma_R^{i+1}(\eta(x)) \\
    %  \sigma_R^i (\zeta(x))
 %\end{pmatrix}
 %=\begin{pmatrix}
      %\sigma_R^i\\
    %  \sigma_R^i  
 %\end{pmatrix}
% \circ
 %\begin{pmatrix}
      %\sigma_R(\eta(x)) \\
    %  \zeta(x)  
 %\end{pmatrix}
 %= \begin{pmatrix}
      %\sigma_R(\eta(y)) \\
    %  \zeta(y)  
 %\end{pmatrix}$ ).
 %Taking all $|R|$ shifts of $x$ (around the range set) into account we get
 %that the number of ways to switch from no discrepancy at site $x$ to a discrepancy
 %at that site is
 %$$[|R|-D(\eta, \zeta)]D(\eta, \sigma_R(\zeta))$$
 %while the number of ways to switch from a discrepancy at $x$ to absence of such is
 %$$D(\eta, \zeta)[|R|-D(\eta, \sigma_R(\zeta))]$$ resulting, after
 %subtraction, in the third term of the RHS in (\ref{SI}). 
 %Due to symmetry, one gets the fourth sum analogously.
 
 Naturally, the first and the second sums on the RHS of (\ref{SI}) represent the contributions made
 to the derivative by all $\begin{pmatrix}
      \sigma \\
      \sigma  
  \end{pmatrix}$ transformations.
 We claim that because $\nu^* \in \mathcal{S}^*$, the first and the second sums on the RHS 
 of (\ref{SI}) must cancel each other. We repeat the same trick: for a range set $R$ containing $x$
 and $\ddot{\eta},\ddot{\zeta} \in \{0,1\}^R$ we consider all shifts 
 $R_{x-\sigma^i(x)}:=R+x-\sigma^i(x)$ of $R$ together with the respective shifts 
 $\begin{pmatrix} \ddot{\eta}^{x-\sigma^i(x)} \\ \ddot{\zeta}^{x-\sigma^i(x)} \end{pmatrix}$ of
 $\begin{pmatrix} \ddot{\eta} \\ \ddot{\zeta} \end{pmatrix}$, for all $i \in \{1,2,\dots,|O(\sigma, x)| \}$.
 For a permutation $\sigma \in \Sigma$ of range $R$, let $\sigma_i$ denote 
 the corresponding $(x-\sigma^i(x))$-shift of $\sigma$. Then $Range(\sigma_i)=R_{x-\sigma^i(x)}$.
 Now, due to the shift-invariant way in which the coupling was constructed,
 $$q^*(\sigma_i, \sigma_i; \ddot{\eta}^{x-\sigma^i(x)}, \ddot{\zeta}^{x-\sigma^i(x)})
 =q^*(\sigma, \sigma; \ddot{\eta}, \ddot{\zeta})$$
 for each $i \in \{1,2,\dots,|O(\sigma, x)| \}$.
 The following are trivial identities: for all $i \in \{1,2,\dots,|O(\sigma, x)| \}$,
 $$\ddot{\eta}^{x-\sigma^i(x)}(x)=\ddot{\eta}(\sigma^i(x)), \qquad
 \ddot{\zeta}^{x-\sigma^i(x)}(x)=\ddot{\zeta}(\sigma^i(x)),$$
 $$\sigma_i(\ddot{\eta})^{x-\sigma^i(x)}(x)=\sigma(\ddot{\eta})(\sigma^i(x))=\ddot{\eta}(\sigma^{i-1}(x))$$ 
 $$\text{ and }
 \sigma_i(\ddot{\zeta})^{x-\sigma^i(x)}(x)=\sigma(\ddot{\zeta})(\sigma^i(x))=\ddot{\zeta}(\sigma^{i-1}(x)).$$
  The total contribution to both first and the second sums on the RHS of (\ref{SI}) made by the transformations
  $\begin{pmatrix}
      \sigma_i \\
      \sigma_i  
  \end{pmatrix}$ for all values of $i \in \{1,2,\dots,|O(\sigma, x)| \}$ is equal to
  $$\sum_{i=1}^{|O(\sigma, x)|} q^*(\sigma_i, \sigma_i; \ddot{\eta}^{x-\sigma^i(x)}, \ddot{\zeta}^{x-\sigma^i(x)})
  \left[ \mathbf{1}_{\{\sigma_i(\ddot{\eta})^{x-\sigma^i(x)}(x) \not= \sigma_i(\ddot{\zeta})^{x-\sigma^i(x)}(x) \}}
  -\mathbf{1}_{\{\ddot{\eta}^{x-\sigma^i(x)}(x) \not= \ddot{\zeta}^{x-\sigma^i(x)}(x)\}} \right]$$
  $$\times \nu^* \left\{ \begin{pmatrix} \eta \\ \zeta \end{pmatrix}: 
  \begin{pmatrix} \eta(R_{x-\sigma^i(x)}) \\ \zeta(R_{x-\sigma^i(x)}) \end{pmatrix} = 
  \begin{pmatrix} \ddot{\eta}^{x-\sigma^i(x)} \\ \ddot{\zeta}^{x-\sigma^i(x)} \end{pmatrix} \right\}$$
  $=q^*(\sigma, \sigma; \ddot{\eta}, \ddot{\zeta}) 
  \cdot \nu^* \left\{ \begin{pmatrix} \eta \\ \zeta \end{pmatrix}:
  \begin{pmatrix} \eta(R) \\ \zeta(R) \end{pmatrix} = \begin{pmatrix} \ddot{\eta} \\ \ddot{\zeta} \end{pmatrix} \right\}$
  $$\times \sum_{i=1}^{|O(\sigma, x)|}
  \left[ \mathbf{1}_{\{ \ddot{\eta}(\sigma^{i-1}(x)) \not= \ddot{\zeta}(\sigma^{i-1}(x)) \}}
  -\mathbf{1}_{\{\ddot{\eta}(\sigma^i(x)) \not= \ddot{\zeta}(\sigma^i(x))\}} \right]=0.$$
  Thus the difference of the first two sums on the RHS of (\ref{SI}) should add up to zero. 
 
 %To prove it, we use the same method as before, where instead of shifting 
 %the permutation $\sigma$, we shift $x$ around the orbit $O(\sigma, x)$, by applying 
 %the permutation over and over again, until returning back to where we started from. So, instead
 %of $x$ we consider all the sites that we have visited,
 %$\sigma(x), \sigma^2(x), \dots, \sigma^{|O(\sigma, x)|}(x)=x$,
 %or, equivalently, we could just consider all the shifts of $\sigma$ so that $x$ is still one of the points of
 %the shifted $O(\sigma, x)$, and shifting $(\eta, \zeta)$ as well. Looking at all the points of the orbit,
 %and a fixed pair $(\eta, \zeta)$,
 %there are as many $i \in \{1,2,\dots,|O(\sigma, x)| \}$ such that 
 %$\{ \eta(\sigma^{i-1}(x)) \not= \zeta(\sigma^{i-1}(x)),\quad \eta(\sigma^i(x)) = \zeta(\sigma^i(x)) \}$
 %as of those satisfying
 %$\{ \eta(\sigma^{i-1}(x)) = \zeta(\sigma^{i-1}(x)),\quad \eta(\sigma^i(x)) \not= \zeta(\sigma^i(x)) \}$ 
 %since we have completed the circle. Hence, after shifting the whole picture $|O(\sigma, x)|$
 %times around the orbit, we have the first two terms canceling each other. 

 Returning to the third and fourth sums on the RHS of (\ref{SI}), since the second sum cancels the first, and
 since the LHS there is $=0$, the third and the fourth sums should also add up to zero.
 We notice that since inside the third sum $\sigma_R(\ddot{\eta}) \geq \ddot{\zeta}$, implying
 $D(\sigma_R(\ddot{\eta}), \ddot{\zeta}) \leq D(\ddot{\eta}, \ddot{\zeta})$, where
 the equality holds only when $D^-(\ddot{\eta},\ddot{\zeta})=0$.
 Similarly $D(\ddot{\eta}, \sigma_R(\ddot{\zeta})) \leq D(\ddot{\eta}, \ddot{\zeta})$
 inside the fourth sum, where
 the equality holds only when $D^+(\ddot{\eta},\ddot{\zeta})=0$.
 That is the number of discrepancies inside $R$ does not change if initially all the discrepancies
 are of the same type, and decreases otherwise. 
 So, $$D(\sigma_R(\ddot{\eta}), \ddot{\zeta}) < D(\ddot{\eta}, \ddot{\zeta}) \text{ in the
 third sum, and } D(\ddot{\eta}, \sigma_R(\ddot{\zeta})) < D(\ddot{\eta}, \ddot{\zeta})$$
 in the fourth sum whenever both types of discrepancies are present inside $R$,
 that is  $D^+(\ddot{\eta},\ddot{\zeta}) \not= 0$ and $D^-(\ddot{\eta},\ddot{\zeta}) \not= 0$.
 Hence for any range set $R$, and any configuration
 $(\ddot{\eta},\ddot{\zeta}) \in \{0,1\}^R \times \{0,1\}^R$ of the coupled process on $R$
 such that $D^+(\ddot{\eta},\ddot{\zeta}) \not= 0$ and $D^-(\ddot{\eta},\ddot{\zeta}) \not= 0$,
 $$\nu^* \Big\{ \begin{pmatrix} \eta \\ \zeta \end{pmatrix}: 
        \begin{pmatrix} \eta(R) \\ \zeta(R) \end{pmatrix} = \begin{pmatrix} \ddot{\eta} \\ \ddot{\zeta} \end{pmatrix}
  \Big\}=0.$$
  Therefore, for all range sets $R$,
   $$\nu^* \Big\{ \begin{pmatrix} \eta \\ \zeta \end{pmatrix}: 
  \quad D^+(\eta(R),\zeta(R)) \not= 0, \quad  D^-(\eta(R),\zeta(R)) \not= 0 \Big\}=0$$
  implying $$\nu^* \big\{(\eta, \zeta)\text{ : } \eta(x)=\zeta(y) \not=\zeta(x)=\eta(y) \big\}=0.$$
  for every $x$ and $y$ in $S$ that both belong to the same range set,\\
  i.e. $\{\sigma \in \Sigma \text{ : } x,y \in Range (\sigma) \} \not= \emptyset$.
  
  The above identity is the first step of the induction. For two points $x$ and $y$ in $S$,
  we let $n(x,y)$ be the least integer $n$ such that there is a sequence
  $$x=x_0, x_1,...,x_n=y$$
  of points in $S$ such that 
  $\{ \sigma \in \Sigma \text{ : } x_{i-1},x_i \in Range(\sigma) \} \not= \emptyset$
  for all $i=1,\dots,n$. Observe that 
  %because of property (\ref{II}), this implies that\\ 
  $\{ \sigma \in \Sigma \text{ : } x_i,x_j \in Range(\sigma) \} = \emptyset$
  for all $0 \leq i,j \leq n$ with $|i-j| \not=1$. We have just proved the basis step $n(x,y)=1$.
  So, for the general step, we assume that the theorem \ref{dplus} is true for $n(x,y)=1,2,\dots,n-1$ (for
  all cases when the connection number $n(x,y)$ that we defined above is any less than 
  the given one). We need to prove that Theorem \ref{dplus} is true for $n(x,y)=n$. We will adapt
  the notation that was used in many papers on interacting particle systems:
  $$\nu^* \left \{ \begin{matrix}
      1\quad  0\\
      0\quad  1\\
      u\quad  v 
      \end{matrix} 
      \right \}
  =\nu^* \big\{(\eta, \zeta)\text{ : } \eta(u)=\zeta(v)=0,\text{ } \zeta(u)=\eta(v)=1 \big\},$$
  for example. Now, for $x$ and $y$ in $S$ with $n(x,y)=n$, we can expand
  \begin{eqnarray*} 
  \nu^* \left \{ \begin{matrix}
      1\quad  0\\
      0\quad  1\\
      x\quad  y 
      \end{matrix} 
      \right \} 
  & = & \nu^* \left\{ \begin{matrix}
      1\quad 1\quad  0\\
      0\quad 1\quad 1\\
      x\quad  x_1\quad y 
      \end{matrix} 
      \right\}
   +\nu^* \left\{ \begin{matrix}
      1\quad 0\quad  0\\
      0\quad 0\quad 1\\
      x\quad  x_1\quad y 
      \end{matrix} 
      \right\} \\
   & + & \nu^* \left\{ \begin{matrix}
      1\quad 1\quad  0\\
      0\quad 0\quad 1\\
      x\quad  x_1\quad y 
      \end{matrix} 
      \right\}
   +\nu^* \left \{ \begin{matrix}
      1 \quad  0 \quad  0\\
      0 \quad  1 \quad 1\\
      x \quad  x_1 \quad y 
      \end{matrix} 
      \right \},
   \end{eqnarray*}
   where the last two terms on the right are equal to zero by the induction hypothesis.
   Here $n(x,x_1)=1$ and $n(x_1,y)=n-1$. Thus, we can show that the first two terms on the RHS are 
   also equal to zero since, by the preceding induction step,
   $$0=\nu^* \left\{ \begin{matrix}
      a_1\quad 1\quad  0\\
      a_2\quad 0\quad 1\\
      x\quad  x_1\quad y 
      \end{matrix} 
      \right\}
      =\nu^* S^*(t) \left\{ \begin{matrix}
      a_1\quad 1\quad  0\\
      a_2\quad 0\quad 1\\
      x\quad  x_1\quad y 
      \end{matrix} 
      \right\}.$$
      Now due to conditions (\ref{I}) and (\ref{II}) there is a $\sigma \in \Sigma$ with 
      $x=x_0,x_1\in Range(\sigma)$ and $x_2,\dots,x_n=y \not\in Range(\sigma)$ such that
      $\sigma(x_0)=x_1$ and $\sigma(x_1)=x_0$ among other things. So,
      $$\nu^* S^*(t) \left\{ \begin{matrix}
      a_1\quad 1\quad  0\\
      a_2\quad 0\quad 1\\
      x\quad  x_1\quad y 
      \end{matrix} 
      \right\}
      \geq \nu^* \left\{ \begin{matrix}
      1\quad a_1\quad  0\\
      0\quad a_2\quad 1\\
      x\quad  x_1\quad y 
      \end{matrix} 
      \right\}  te^{-ct}q(\sigma),$$
      where the constant $c$ is greater than the sum of the rates of all other permutations in $\Sigma$
      containing any of the $x_i$'s in their ranges.

      Observe that one does not really need
      $\sigma(x_1)=x_0$ when doing this proof with weaker conditions than (\ref{I}) and (\ref{II})
      that were mentioned in \ref{conditions}.
      
      So,
      $$  \nu^* \left\{ \begin{matrix}
      1\quad  0\\
      0\quad  1\\
      x\quad  y 
      \end{matrix} 
      \right\} =0$$
      for all $x$ and $y$ in $S$ with all values of $n(x,y)$, and Theorem \ref{dplus} is proved. $\square$
      
\subsection{Proof of Theorem \ref{ISe} }
   Since Theorem \ref{star} and Theorem \ref{dplus} are now proved, 
   the proof of Theorem \ref{ISe} is word to word
   identical to the analogous case in the theory of exclusion processes and is a part of the system of
   results developed by T.Liggett for the exclusion processes that we are trying to redo for the permutation 
   processes. Though since the proof is short, and since we need to inform the reader of
   why Theorem \ref{star} and Theorem \ref{dplus} are so important as parts of the proof of Theorem \ref{ISe},
   we are going to copy the proof in the remaining few lines of this section.
   
   {\it Proof of Theorem \ref{ISe}: }
   Since $\int \Omega f d\nu_{\rho}=0$ , $\nu_{\rho} \in \mathcal{I}$ and obviously
   $\nu_{\rho} \in \mathcal{S}$ for all $0 \leq \rho \leq 1$. Furthermore, 
   $\nu_{\rho} \in \mathcal{S}_e$, since it is spatially ergodic. Therefore,
   $\nu_{\rho} \in (\mathcal{I} \cap \mathcal{S})_e$.
   
   For the converse, take $\nu \in (\mathcal{I} \cap \mathcal{S})_e$. By Theorem \ref{star}(e),
   for any $0 \leq \rho \leq 1$, there is a $\nu^* \in (\mathcal{I}^* \cap \mathcal{S}^*)_e$
   with marginals $\nu_{\rho}$ and $\nu$. By Theorem \ref{dplus},
   $$\nu^* \big\{(\eta, \zeta)\text{ : } \eta \leq \zeta \quad \eta \not= \zeta \big\}
   +\nu^* \big\{(\eta, \zeta)\text{ : } \zeta \leq \eta \quad \eta \not= \zeta \big\}
   +\nu^* \big\{(\eta, \zeta)\text{ : } \eta = \zeta \big\} =1.$$
   Since the three sets above are closed for the evolution and translation invariant, and
   since $\nu^*$ is extremal, it follows that one of the three sets has full measure.
   Therefore, for every $0 \leq \rho \leq 1$, either $\nu \leq \nu_{\rho}$ or
   $\nu_{\rho} \leq \nu$. It follows that $\nu = \nu_{\rho_0}$ where $\rho_0$ is determined
   by
   $$\nu \leq \nu_{\rho} \quad \text{for } \rho > \rho_0,$$
   $$\nu \geq \nu_{\rho} \quad \text{for } \rho < \rho_0.$$
   \qquad\qquad\qquad\qquad\qquad\qquad\qquad\qquad\qquad\qquad\qquad\qquad\qquad\qquad $\square$

\section*{Acknowledgments}
 The author wishes to thank Tom Liggett who suggested that the author investigates the model,
 and who was the main source of help and inspiration along the way.
    
\bibliographystyle{amsplain}

 \textit{Yevgeniy Kovchegov\\
 Department of Mathematics, UCLA\\
 Email: yevgeniy@math.ucla.edu\\
 Fax: 1-310-206-6673}

 \end{document}